\documentclass[times,10pt]{article}
\usepackage{amssymb,amsmath}

\begin{document}
\centerline{\bf A Note on the Riemann $\xi-$Function}\vskip .4in

\centerline{M.L. GLasser}\vskip .1in

\centerline{Department of Theoretical Physics, University of Valladolid, Valladolid (Spain)}\vskip .2in
\centerline{Department of Physics, Clarkson University, Potsdam, NY(USA)}\vskip 1in

\centerline{\bf Abstract}
\vskip .1in

This  note investigates a number of integrals of and integral equations satisfied by Riemann's $\xi-$function. A different,  less restrictive, derivation of one of his key identities is provided. This work centers on the critical strip and it is argued that the line  $s=3/2+i t$ , e.g., contains a holographic image of this region. \vskip 1in

\noindent
MSC classes: 11M06, 11M26, 11M99, 26A09, 30B40, 30E20, 30C15, 33C47, 33B99, 33F99
\vskip .8in
\noindent
{\bf Keywords:} $\xi$ function, $\zeta$. function, Integral equations

\vskip 1in
\noindent
email: laryg@clarkson.edu 

\newpage

\centerline{\bf Introduction}\vskip .2in

The notation used throughout this note is:
$$\xi(s)=(s-1)\pi^{-s/2}\Gamma(1+s/2)\zeta(s)$$
$$\rho=\sigma+i\tau$$
$$\Xi(\tau)=\xi(1/2+i\tau)$$
$$E_z(a)=\int_1^{\infty}\frac{dt}{t^z}e^{-at},$$
$$\psi(x)=\sum_{n=1}^{\infty}e^{-\pi n^2 x}=\frac{1}{2}[\theta_3(0,e^{-\pi x})-1]$$
$$ \quad J(\rho)=\int_0^1 dt [t^{\rho-2}+t^{(1-\rho)-2}]\psi(1/t^2).$$
$\gamma$ is the contour consisting of the two parallel lines $[c-i\infty,c+i\infty]$, $[1-c+i\infty,1-c-i\infty]$, $1<c<2$, which span the {\it critical strip} $0<\sigma<1$, $\rho_n=1/2+i\tau_n$ is the n-th zero of $\zeta(s)$ on the critical line in the upper half plane.\vskip .2in

The function $\xi(s)$, introduced by Riemann[1], satisfies the simple functional equation $\xi(1-s)=\xi(s)$, is analytic, decays exponentially as $|\rho|\rightarrow\infty$ and possesses the same zeros  in the critical strip as $\zeta(s)$.  By Cauchy's theorem, one has, for $0<Re[s]<1$
$$\xi(s)=\int_{\gamma}\frac{dt}{2\pi i}\frac{\xi(t)}{t-s},\eqno(1)$$
which, in view of the functional equation, can be written
$$\xi(s)=\int_{c-i\infty}^{c+i\infty}\frac{dt}{2\pi i}\xi(t)\left(\frac{1}{t-s}+\frac{1}{t-1+s}\right)\eqno(2)$$
and expresses the values of $\xi$  inside the critical strip entirely in terms of its values in a region where $\zeta(s)$ is completely known from its defining series, say. In the following section this feature will be exploited to obtain several known and some, perhaps, unfamiliar identities.\vskip .2in

\centerline{\bf Calculation}\vskip .1in

We begin by recalling the tabulated inverse Mellin  transform[2]

$$\int_{c-i\infty}^{c+i\infty}\frac{dt}{2\pi i}x^{-t}\Gamma(t)\zeta(2t)=\sum_{n=1}^{\infty}e^{-n^2 x}\eqno(3)$$
from which, by differentiation, one finds the useful inverse Mellin transform
$$F(x)=\int_{c-i\infty}^{c+i\infty}\frac{dt}{2\pi i}x^{-t}\xi(t)=4\pi^2 x^4\sum_{n=1}^{\infty}n^4 e^{-\pi n^2 x^2}-6\pi x^2\sum_{n=1}^{\infty}n^2e^{-\pi. n^2x^2},\eqno(4)$$
Eq.(4) has  been presented, in different form, by Patkowski[3].

Parenthetically, we note that if $f$ is integrable and odd, then $f(1-2t)\xi(t)$ is odd under $t\rightarrow(1-t)$ so that

$$\int_{c-i\infty}^{c+i\infty}\frac{dt}{2\pi i}f(1-2t)\xi(t)=\int_{\gamma}\frac{dt}{2\pi i}f(1-2t)\xi(t)=0.\eqno(5)$$
Thus, by noting Romik's formulas[4] for the values of the Theta function $\theta_3(0,q)$ and its derivatives we have, from(4)
 $$ \int_{c-i\infty}^{c+i\infty}\frac{dt}{2\pi i}\frac{\xi(t)}{t}=\frac{1}{2}-\frac{\Gamma(5/4)}{\sqrt{2}\pi^{3/4}}\eqno(6).$$
and from (5)

$$\int_{c-i\infty}^{c+i\infty}\frac{dt}{2\pi i}\xi(t)\left[\frac{(1-2t)}{4\tau_n^2+(1-2t)^2}\right]=0\eqno(7)$$
 $$\int_{c-i\infty}^{c+i\infty}\frac{dt}{2\pi i}\xi(t)(2t-1)^{2n+1}=0,\quad n=0,1,2,\cdots..\eqno(8a)$$
  $$\int_{c-i\infty}^{c+i\infty}\frac{dt}{2\pi i}t\xi(t)=\frac{1}{2}\int_{c-i\infty}^{c+i\infty}\frac{dt}{2\pi i}\xi(t)=\eqno(8b)$$
    $$\pi^2\sum_{n=-\infty}^{\infty} n^4e^{-\pi n^2}-\frac{3\pi}{2}\sum_{n=-\infty}^{\infty}n^2e^{-\pi n^2}\eqno(8c)$$
$$=\frac{\Gamma(5/4)}{128\sqrt{2}\pi^{19/4}}[\Gamma^8(1/4)-96\pi^4].\eqno(8d)$$
 $$\int_{c-i\infty}^{c+i\infty}\frac{dt}{2\pi i}\frac{\xi(t)}{t(1-t)}=\frac{1}{2}\left(1-\frac{\pi^{1/4}}{\Gamma(3/4)}\right)=\int_{1-c-i\infty}^{1-c+i\infty}\frac{dt}{2\pi i}\frac{\xi(t)}{t(1-t)}\eqno(9)$$
 None of these appears to have been recorded previously.
 
 Next, by rewriting (2), we have \vskip .2in
 
 \noindent
 {\bf Theorem 1}\vskip .1in
 Within the critical strip Riemann's function $\xi(s)$ obeys the integral equation

  $$\xi(s)=1-\frac{\pi^{1/4}}{2\Gamma[3/4)}-\int_{c-i\infty}^{c+i\infty}\frac{dt}{2\pi i}\frac{\xi(t)}{t}\left[\frac{2s(1-s)-t}{s(1-s)-t(1-t)}\right], \quad 1<c<2.\eqno(10)$$
  or
  
  $$\xi(s)=\frac{1}{2}+\int_{c-i\infty}^{c+i\infty}\frac{dt}{2\pi i}\frac{\xi(t)}{t}\left[\frac{1}{1-t}-\frac{2s(1-s)-t}{s(1-s)-t(1-t)}\right].\eqno(11)$$\vskip .2in
   
\noindent
From this, one finds \vskip .1in

\noindent
{\bf Corollary 1}\vskip .1in
       
   $$\xi(s)=2\pi^2\sum_{n=1}^{\infty}\int_1^{\infty}dt\left(t^{s/2}+t^{(1-s)/2}\right)\left(n^4t-\frac{3}{2\pi}n^2\right)e^{-n^2\pi t}.\eqno(12)$$
   
   $$=\frac{\pi^{1/4}}{2\Gamma(3/4)}-
   \pi\sum_{n=1}^{\infty}n^2\left[sE_{(1-s)/2}(\pi n^2)+(1-s)E_{-s/2}(\pi n^2)\right],\eqno(13)$$
   \vskip .1in
   
   Eq. (10) is equivalent to the very important eq(3.10) in Milgram's paper[5] and  (13), apart from having summed a series, is LeClair's key formula (15) in [6]\vskip .2in
   
   To explore further consequences of (2), note that
   the Mellin transform
   $$\phi(x)= \int_{c-i\infty}^{c+i\infty}\frac{dt}{2\pi i}x^{-t}\frac{\xi(t)}{t-s}\eqno(14)$$
   satisfies the linear differential equation
   $$\phi'(x)+\frac{s}{x}\phi(x)=-\frac{1}{x}F(x),\quad \phi(\infty)=0\eqno(15)$$
   where $F$ is defined in (4), so after a bit of easy analysis
   
   $$\phi(x)=2\pi^2\sum_{n=1}^{\infty} E_{\frac{x}{2}+1}(\pi n^2)-3\pi\sum_{n=1}^{\infty}E_{\frac{x}{2}}(\pi n^2).\eqno(16)$$

    By applying (16)  to (2) one has (note that $\tau$ here is not restricted to be real)\vskip .2in
    
    \noindent{\bf Theorem 2}
    \vskip .1in
   In the critical strip , i.e. for $|Im[\tau]|<1/2$, Riemann's function $\Xi(\tau)$ satisfies the integral equation
    
$$\Xi(\tau)=\frac{1}{\pi i}\int_{-\infty+ic}^{\infty+ic}\frac{t\,  \Xi(t)}{t^2-\tau^2}dt,\quad -3/2<c<-1/2.\eqno(17)$$\vskip .2in

\noindent
{\bf Corollary 2}\vskip .1in

   $$\Xi(\tau)=4\pi^2\sum_{n=1}^{\infty}\int_1^{\infty}dt\, t^{1/4}\cos(\tau \ln\sqrt{t})\left(n^4t-\frac{3}{2\pi}n^2\right)e^{-n^2\pi t}.\eqno(17a)$$
   $$=4\int_1^{\infty}dt\, t^{1/4}\cos[\frac{1}{2}\tau\ln t][t\psi''(t)+\frac{3}{2}\psi'(t)]\eqno(17b)$$
   
   $$=\frac{1}{2}-(\tau^2+1/4)\sum_{n=1}^{\infty}Re\,  E_{\frac{3}{4}+i\frac{\tau}{2}} (\pi n^2)\eqno(17c)$$
     So
   
   $$\Xi(\tau)=1/2-(\tau^2+1/4)\int_1^{\infty}\frac{dt}{t^{3/4}}\cos(\frac{\tau}{2}\ln t)\psi(t)\eqno(18)$$
    \vskip .2in
    
    Now, (18), which appears in Riemann's paper[7], can be rewritten
    $$\xi(\rho)=\frac{1}{2}-(\alpha+i\beta)\int_0^1dt [t^{\rho-2}+t^{(1-\rho)-2}]\psi(1/t^2).\eqno(19)$$
    where $\alpha=\sigma(1-\sigma)+\tau^2$ and $\beta=(1-2\sigma)\tau.$
   For $\sigma=1/2$, thus restricting $\tau$ to real values, (18) gives\vskip .1in
   
   \noindent
   {\bf Corollary 3}\vskip .1in
   $\rho$ is a zero of $\zeta(s)$
    on the critical line, $\sigma=1/2$, if and only if,
   $$Re\int_0^{1}t^{\rho-2}\psi(1/t^2)dt=\frac{1}{4|\rho|^2}.\eqno(20)$$
   \vskip .1in
   In our key result (17) for simplicity we choose $c=1$ and to emphasize that $\tau$ is not restricted to real values, write 
   $$\xi(\rho)=\frac{1}{\pi }\int_{-\infty}^{\infty}dt \xi(3/2+i t)\frac{1+it}{(1+it)^2+[\frac{1}{2}-\sigma+i\tau)^2}\eqno(21)$$
   \vskip .2in
   
   Before leaving this section, we note that (17) is the source of a great number of fascinating integral identities found by multiplying both sides by a suitable function $g(\tau)$  and integrating over $\tau$. In the next section a few of these are presented omitting details.\vskip .2in
   
   \centerline{\bf Additional integrals}\vskip .1in
   
  $$ \int_0^{\infty}\cos(xt)\Xi(t)dt=\frac{1}{2}e^{-x}\int_{-\infty}^{\infty}e^{-ixt}\xi(3/2+it)dt,\quad x>0\eqno(i),$$
  so, from [2], for $x>0$
  $$\int_{-\infty}^{\infty}e^{-ixt}\xi(3/2+it)dt=8\pi e^{-3x/2}[e^{-2x}\psi''(e^{-2x})-\frac{3}{2}\psi'(e^{-2x})]\eqno(ii)$$
  and
  $$\int_{-\infty}^{\infty}[\xi(1/2+i t)-\xi(3/2+it)]dt=0.\eqno(iii)$$
  Similarly,
  $$\int_0^{\infty}\Xi(t)\frac{\cos(xt)}{t^2+b^2}dt=\frac{e^{-x}}{2b}\int_{-\infty}^{\infty}dt\frac{\xi(3/2+it)}{(1+it)^2+b^2}[(1+it)e^{(1-b)x}-be^{-ixt}],\eqno(iv)$$
   so by [9]
  $$\int_{-\infty}^{\infty}dt\frac{\xi(3/2+it)}{(1+it)^2+1/4}[(1+it)e^{x/2}-\frac{1}{2}e^{-ixt}]=\frac{\pi}{2}[e^{x/2}[e^x-2\psi(e^{-2x})].\eqno(v)$$
  Thus
  $$\int_{-\infty}^{\infty}dt\frac{\xi(3/2+it)}{(1+it)^2+1/4}(\frac{1}{2}+it)=\frac{\pi}{2}(1-\sqrt{2}\Gamma(1/4)\pi^{-3/4}).\eqno(vi)$$
  $$\int_0^{\infty} J_0(a t)\Xi(t)dt=\frac{e^{-x}}{2}\int_{-\infty}^{\infty}\xi(3/2+it)I_0[a(1+it)]dt.\eqno(vii)$$
  $$\int_0^T\Xi(t)dt=\frac{1}{\pi}\int_{-\infty}^{\infty}\xi(3/2+it)\tan^{-1}\left(\frac{T}{1+it}\right).\eqno(viii)$$
  Eq.(viii) is equivalent to
  $$\int_0^T\xi(\sigma+it)dt=\frac{1}{\pi}\int_{-\infty}^{\infty}dt\xi(3/2+it)\tan^{-1}\left(\frac{T+
 i( \frac{1}{2}-\sigma)}{1+it}\right)+g_0(\sigma)\eqno(ix)$$
  where $0<\sigma<1$ and
  $$g_0(\sigma)=(\sigma-1/2)i\int_0^1\xi(1/2+(\sigma-1/2) t]dt\eqno(ixa)$$
  is independent of $T$.
  What  these examples have in common is that all information on the critical line is equivalent to information on the line $\sigma=3/2$.
\vskip .2in
   
   \centerline{\bf Discussion}\vskip .1in
   
   The thrust of the preceding sections, and the key formulas  have all been confirmed numerically with Mathematica,
 is that all the features of $\xi$ in the critical strip are encoded on the lines $\sigma=c$, $1<c<2$ in a holographic manner made explicit by (21). This has the form
   $$\xi(\rho)=\int \xi(c+it)R(\sigma,\tau;t)dt\eqno(22)$$
   where R is a {\it rational} function. As functions go, $\xi(3/2+it)$ lis uncomplicated: it has no zeros or poles,  it is infinitely differentiable and it  decays exponentially. Its real and imaginary parts are even and odd, respectively, and related  to each other by the Cauchy-Riemann equations. Both of the latter functions are, except near $t=0$, almost featureless. Furthermore, the zeta function component is nearly equal to unity over most of the integration range. Thus any interesting feature in the critical strip must be ascribed largely to features of $R$ which  ought be trivial to analyze.   
        
       Equation  (20) (in essence due to Riemann)  has been confirmed  for the first 1000 critical zeros, which are  available on Mathematica, in its exact form, but even if $\psi(x)$ is truncated to one exponential it  is satisfied to many decimal place accuracy for large magnitude zeros. In this case, the integral is  $E_{3/4+i\tau_n/2}(\pi)$ and  by asymptotic expansion should be capable of producing a formula for $\tau_n$ similar to LeClair's[6] and Milgram's[7], but with less complexity.     That is, to derive an expression for $\rho_n$,  (20)  one expresses in the form
        $$Re\left[e^{-(\rho+1/2)\ln\sqrt{\pi}}\Gamma\left[\frac{1}{2}(\rho+\frac{1}{2})\right]+g(\rho)\right]=0\eqno(23 a)$$
       $$g(\rho)=\frac{\rho+1/2}{|\rho+1/2|^2}\;_1F_1(\rho+1/2;\rho+3/2;-\pi)-\frac{1}{4|\rho|^2}\eqno(23b)$$
 Now, along the critical line   the real part of the function $g( \rho)$ in (23b) is non-oscillatory, monotonically decreasing and smaller than the accuracy we are trying to achieve, although much larger than the first 	term  in (23a).which is oscillatory. However, in the spirit of [6 ] one ignores $g$ and thus approximates (23a) as
 $$Re\left[e^{-i\tau\ln(\sqrt{\pi})}\Gamma\left(\frac{1}{2}+i\frac{\tau}{2}\right)\right]=0\eqno(23c)$$
 Following  [6] one now applies Stirling's formula in (22) and solves for $\tau_n$ to obtain an analogue of LeClair's formula[6. (20)]. 
 
 Finally, as anyone, who writes on the zeta function  must find irresistible, I conclude this note  with speculations on the elephant in the room, the Riemann hypothesis(RH). In essence, the latter claims nothing more than that $\xi(\rho)$ does not vanish for $0<\sigma<1/2$. However, as remarked above, this is  a question of the simultaneous vanishing of the real and imaginary parts of the integral (22) for which the zeta function itself is nearly irrelevant. For $\sigma=1/2$, the imaginary part goes away and it is known that the real part has countably many roots $\tau_n$. Otherwise, the situation appears to depend mainly on the nature of the rational function $R$, which depends on $\sigma$, than  on $\xi$, which does not. so the validity of the RH  should be easy to resolve. It could be that the vast literature concerning $\zeta$ in the critical strip is a red herring.
      
        Another possible approach is to break the integral  $J(\rho)$ in (19)  into real and imaginary parts. Then if  $\rho$ is a zero of $\zeta$ with $0<\sigma<1/2$,
        $$\alpha Re[J]-\beta Im[J]=\frac{1}{2}\eqno(24a)$$
        $$\beta Re[J]+\alpha Im[J]=0\eqno(23b)$$
These two relations lead to
\vskip .1in
\noindent
{\bf Corollary 4}\vskip .1in
The Riemann conjecture is true if and only if
$$Re\int_0^1 dt [t^{\rho-2}+t^{(1-\rho)-2}]\psi(1/t^2)=\frac{\alpha}{2(\alpha^2-\beta^2)}\eqno(25)$$
has no solution for $0<\sigma<1/2$.\vskip .1in

 Since the critical strip is known to be free of such zeros to astronomical values of $|\rho|$, the resolution of this matter might  be settled by extracting a low order asymptotic estimate of the Mellin transform $J(\rho)$ and analyzing the resulting algebraic equation. Since truncating $\psi$ to a single term appears to yield accurate results for large $n$, it  may be reasonable to  conjecture:\vskip .2in
 The Riemann hypothesis is true if for large $t$ the equation
$$\frac{(1-s) s+t^2}{\left((1-s) s+t^2\right)^2-(1-2 s)^2
   t^2}-\Re\left(E_{-\frac{s}{2}-\frac{i t}{2}+1}(\pi )+E_{\frac{1}{2} (s+i t+1)}(\pi
   )\right)=0$$
   has no real solution $t$ for $0<s<1/2$. However, such simple expedients tend to be illusory since no matter how small it is, a positive number is not zero.
      \vskip .1in

  \newpage      
         
 \vskip .4in\noindent{\bf Acknowledgements}\vskip ,1in
 
 This work was inspired by reading [5] and I thank its author for  this opportunity. I also thank Dr. Michael Milgram and Prof. Richard Brent for valuable comments and suggestions. This work was partially supported by the Spanish  grant MINECO (MTM2014-57129-C2-1-P), Junta de Castilla y Le\'on, FEDER projects (BU229P18, VA057U16, and VA137G18).\vskip .2in

 \centerline{\bf References}\vskip .1in
 
 \noindent
 [1] E.C. Titchmarsh and D.R Heath-Brown. The Theory of the Riemann Zeta-Function. Oxford
Science Publications, Oxford, Second edition, 1986. \vskip .1in
 
 \noindent
 [2] Ibid (2.16.1)\vskip .1in
 
 \noindent
 [3] A.E. Patkowski, A New Integral Equation and some Integrals Assviated with Number Theory, arXiv e-prints, March 2018. arXiv: 1407.2983v6\vskip .1in
 
 \noindent
[4] Dan Romik. The Taylor coefficients of the Jacobi theta constant  $\theta_3$, arXiv e-prints, July
2018. arXiv:1807.06130.\vskip .1in
 
 \noindent
 [5] Michael Milgram, An Integral Equation for Riemann?s Zeta Function and its Approximate Solution,
      arXiv e prints, January 20191901.01256v1 [math.CA]\vskip .1in
      
      \noindent
 [6]  Andre LeClair. An electrostatic depiction of the validity of the Riemann Hypothesis and a
formula for the N-th zero at large N. Int. J. Mod. Phys., A28:1350151, 2013. Also available
from http://arxiv.org/abs/1305.2613v3\vskip .1in
 
 \noindent
 [7] Michael Milgram,  Exploring Riemann?s functional equation, Cogent Math. (2016), 3: 1179246
      http://dx.doi.org/10.1080/23311835.2016.1179246\vskip .1in
      
 \noindent
 [8]  H.M. Edwards. Riemann's Zeta Function. Dover, 2001.   \vskip .1in
 
 \noindent
 [9] (2.16.2) of reference[1].  \vskip .1in

\end{document}